\documentclass{article}

\newcommand{\qed}{\hfill $\dashv$}
\newcommand{\cirk}{\,{\raisebox{.3ex}{\tiny $\circ$}}\,}

\begin{document}

\title{Coherent Bicartesian and Sesquicartesian Categories}
\author{{\sc Kosta Do\v sen} and {\sc Zoran Petri\' c}
\\[1ex]
{\small Mathematical Institute, SANU}\\[-.5ex]
{\small Knez Mihailova 35, p.f. 367, 11001 Belgrade,
Serbia}\\[-.5ex]
{\small email: \{kosta, zpetric\}@mi.sanu.ac.yu}}
\date{\small May 2006}
\maketitle

\begin{abstract}
\noindent Coherence is here demonstrated for sesquicartesian
categories, which are categories with nonempty finite products and
arbitrary finite sums, including the empty sum, where moreover the
first and the second projection from the product of the initial
object with itself are the same. (Every bicartesian closed
category, and, in particular, the category {\bf Set}, is such a
category.) This coherence amounts to the existence of a faithful
functor from categories of this sort freely generated by sets of
objects to the category of relations on finite ordinals, and it
yields a very easy decision procedure for equality of arrows.
Restricted coherence holds also for bicartesian categories where,
in addition to this equality for projections, we have that the
first and the second injection to the sum of the terminal object
with itself are the same.

The printed version of this paper (in: R.\ Kahle et al.\ eds, {\it
Proof Theory in Computer Science}, Lecture Notes in Computer
Science, vol.\ 2183, Springer, Berlin, 2001, pp.\ 78-92) and
versions previously posted here purported to prove unrestricted
coherence for the bicartesian categories mentioned above. Lemma
5.1 of these versions, on which the proof of coherence for
sesquicartesian categories relied too, is however not correct. The
present version of the paper differs from the previous ones also
in terminology.
\end{abstract}

\vspace{0.3cm}
\noindent {\footnotesize Mathematics Subject Classification (2000): 18A30,
18A15, 03G30, 03F05} \vspace{0.5cm}

\section{Introduction}

\noindent The connectives of conjunction and disjunction in classical and
intuitionistic logic make a structure corresponding to a distributive
lattice. Among nonclassical logics, and, in particular, among substructural
logics, one finds also conjunction and disjunction that make a structure
corresponding to a free lattice, which is not distributive.

In proof theory, however, we are not concerned only with a consequence {\em %
relation}, which in the case of conjunction and disjunction gives rise to
a lattice
order, but we may distinguish certain deductions with the same premise and
the same conclusion. For example, there are two different deductions from $%
A\wedge A$ to $A$, one corresponding to the first projection and the other
to the second projection, and two different deductions from $A$ to $A\vee A$,
one corresponding to the first injection and the other to the second
injection.

If we identify deductions guided by normalization, or cut elimination, we
will indeed distinguish the members in each of these two pairs of
deductions, and we will end up with natural sorts of categories, whose
arrows stand for deductions. Instead of nondistributive lattices, we obtain
then categories with binary products and sums (i.e. coproducts), where
the product $\times$ corresponds to conjunction and the sum $+$ to
disjunction. If, in
order to have all finite products and sums, we add also the empty product,
i.e. the terminal object, which corresponds to the constant true
proposition, and the empty sum, i.e. the initial object, which corresponds
to the constant absurd proposition, we obtain {\em bicartesian} categories,
namely categories that are at the same time {\em cartesian}, namely, with
all finite products, and {\em cocartesian}, with all finite sums.
Bicartesian categories need not be distributive.

One may then enquire how useful is the representation of deductions
involving conjunction and disjunction in bicartesian categories for
determining which deductions are equal and which are not. A drawback is that
here cut-free, i.e. composition-free, form is not unique. For example, for $%
f:A\rightarrow C$, $g:A\rightarrow D$, $h:B\rightarrow C$ and
$j:B\rightarrow D$ we have the
following two composition-free arrow terms
\[
\frac{\langle f,g\rangle :A\rightarrow C\times D\quad \quad \langle
h,j\rangle :B\rightarrow C\times D}{[\langle f,g\rangle ,\langle h,j\rangle
]:A + B\rightarrow C\times D}
\]
\[
\frac{\lbrack f,h]:A + B\rightarrow C\quad \quad [g,j]:A + B\rightarrow D%
}{\langle [f,h],[g,j]\rangle :A + B\rightarrow C\times D}
\]
which designate the same arrow, but it is not clear which one of them
is to be
considered in normal form. The same problem arises also for distributive
bicartesian categories, i.e. for classical and intuitionistic conjunction
and disjunction. (This problem is related to questions that arise in natural
deduction with Prawitz's reductions tied to disjunction elimination,
sometimes called ``commuting conversions''.)

Cut elimination may then be supplemented with a {\em coherence}
result. The notion of coherence is understood here in the
following sense. We say that coherence holds for some sort of
category iff for a category of this sort freely generated by a set
of objects there is a faithful functor from it to a {\em graphical
}category, whose arrows are sets of links between the generating
objects. These links can be drawn, and that is why we call the
target category ``graphical''. Intuitively, these links connect
the objects that must remain the same after ``generalizing'' the
arrows. We shall define precisely below the graphical category we
shall use for our coherence results, and the intuitive notion of
generalizing (which originates in \cite {lam68} and \cite{lam69})
will become clear. This category will be the category of relations
on finite ordinals. In other cases, however, we might have a
different, but similar, category. It is desirable that the
graphical category be such that through coherence we obtain a
decision procedure for equality of arrows, and in good cases, such
as those investigated in \cite {dos00b} and here, this happens
indeed.

Although this understanding of coherence need not be the most standard one,
the paradigmatic results on coherence of \cite{mac63} and \cite{kel71} can
be understood as faithfulness results of the type mentioned above, and this
is how we understand coherence here. We refer to \cite{dos00b}, and papers
cited therein, for further background motivation on coherence.

This paper is a companion to \cite{dos00b}, where it was shown that
coherence holds for categories with binary, i.e. nonempty finite, products
and sums, but without the terminal and the initial object, and without
distribution. One obtains thereby a very easy decision procedure for
equality of arrows. With the help of this coherence, it was also
demonstrated that the categories in question are maximal, in the sense that
in any such category that is not a preorder all the equations between arrows
inherited from the free categories with binary products and sums are the
same. An analogous maximality can be established for cartesian categories
(see \cite{dos00}) and cartesian closed categories (see \cite{sim95} and
\cite{dos00a}).

In this paper we shall be concerned with categories with nonempty finite
products and arbitrary finite sums, including the empty sum, i.e. the
initial object. We call such categories {\em sesquicartesian}. As a matter
of fact, {\em sesquicocartesian} would be a more appropriate label for these
categories, since they are {\em cocartesian} categories---namely, categories
with arbitrary finite sums---to which we add a part of the cartesian
structure---namely, nonempty finite products. {\em Sesquicartesian} is more
appropriate as a label for cartesian categories---namely, categories with
arbitrary finite products---to which we add nonempty finite sums, which are
a part of the cocartesian structure. However, sesquicocartesian and
sesquicartesian categories so understood are dual to each other, and in a
context where both types of categories are not considered, there is no
necessity to burden oneself with the distinction, and make a strange new
name even stranger. So we call here {\em sesquicartesian} categories with
nonempty finite products and arbitrary finite sums.

We shall show that coherence holds for sesquicartesian categories
in which the first and the second projection arrow from the
product of the initial object with itself are equal. Such
sesquicartesian categories were called \emph{coherent}
sesquicartesian categories in the printed version of this paper.
Now we use just \emph{sesquicartesian category} to designate what
we used to call \emph{coherent sesquicartesian category}. Every
bicartesian closed category, and, in particular, the
sesquicartesian category {\bf Set} of sets with functions as
arrows, is a sesquicartesian category in this new sense of the
term. It is not true, however, that sesquicartesian categories are
maximal, in the sense in which cartesian categories, cartesian
closed categories and categories with binary products and sums are
maximal.

Bicartesian categories are also not maximal. Coherence does not
hold for bicartesian categories in general. We will prove however
a restricted coherence result for bicartesian categories where
besides the equality of the projection arrows mentioned above we
also have that the first and the second injection arrow to the sum
of the terminal object with itself are equal. We call such
bicartesian categories \emph{dicartesian} categories. (In the
printed version of this paper we called them \emph{coherent}
bicartesian categories.) The bicartesian category {\bf Set} is not
such a category, but a few natural examples of such categories may
be found in Section 7. Such categories are also not maximal.

As in \cite{dos00b}, our coherence for sesquicartesian categories
yields a very easy decision procedure for equality of arrows in
the categories of this kind freely generated by sets of objects.
For dicartesian categories this procedure is partial. Without
maximality, however, the application of this decision procedure to
an arbitrary category of the appropriate kind is limited. We can
use this decision procedure only to show that two arrows inherited
from the free category are equal, and we cannot use it to show
that they are not equal.

We said that this paper is a companion to \cite{dos00b}, but except for
further motivation, for which we refer to \cite{dos00b}, we shall strive to
make the present paper self-contained. So we have included here definitions
and proofs that are just versions of material that may be found also in \cite
{dos00b}, but which for the sake of clarity it is better to adapt to the new
context.

\section{Free Dicartesian and Ses\-qui\-cartesian Catego\-ries}

The propositional language ${\cal P}$ is generated from a set of {\em %
propositional letters} ${\cal L}$ with the nullary connectives, i.e.
propositional constants, I and O, and the binary connectives $\times$ and $+$%
. The fragment ${\cal P}_{\times,+,\mbox{\scriptsize{\rm O}}}$ of ${\cal P}$
is obtained by omitting all formulae that contain I. For the propositional
letters of ${\cal P}$, i.e. for the members of ${\cal L}$, we use the
schematic letters $p,q,\ldots$, and for the formulae of ${\cal P}$, or of
its fragments, we use the schematic letters $A,B,\ldots,A_1,\ldots$ The
propositional letters and the constants I and O are {\em atomic} formulae.
The formulae of ${\cal P}$ in which no propositional letter occurs will be
called {\em constant objects}.

Next we define inductively the {\em terms} that will stand for the arrows of
the free {\em dicartesian category} ${\cal D}$\ generated by ${\cal %
L}$. Every term has a {\em type}, which is a pair $(A,B)$ of formulae of $%
{\cal P}$. That a term $f$ is of type $(A,B)$ is written $f:A\rightarrow B$.
The {\em atomic} terms are for every $A$ and every $B$ of ${\cal P}$
\[
\begin{array}{lcl}
& \mbox{\bf 1}_A:A\rightarrow A, &  \\
k_A:A\rightarrow\mbox{\rm I}, &  & l_A:\mbox{\rm O}\rightarrow A, \\
k^1_{A,B}:A\times B\rightarrow A, &  & l^1_{A,B}:A\rightarrow A+B, \\
k^2_{A,B}:A\times B\rightarrow B, &  & l^2_{A,B}:B\rightarrow A+B, \\
w_A:A\rightarrow A\times A, &  & m_A:A+A\rightarrow A.
\end{array}
\]
The terms $\mbox{\bf 1}_A$ are called {\em identities}. The other terms of $%
{\cal D}$\ are generated with the following operations on terms,
which we present by rules so that from the terms in the premises
we obtain the terms in the conclusion:

\[
\frac{f:A\rightarrow B\quad \quad g:B\rightarrow C}{g\cirk
f:A\rightarrow C}
\]
\[
\frac{f:A\rightarrow B\quad \quad g:C\rightarrow D}{f\times g:A\times
C\rightarrow B\times D}\quad \quad \quad \frac{f:A\rightarrow B\quad \quad
g:C\rightarrow D}{f+g:A+C\rightarrow B+D}
\]

\noindent We use $f,g,\ldots,f_1,\ldots$ as schematic letters for terms of $%
{\cal D}$.

The category ${\cal D}$\ has as objects the formulae of ${\cal P}$
and as arrows equivalence classes of terms so that the following
equations are satisfied for $i\in\{1,2\}$:
\[
\begin{array}{l}
(cat\: 1)\quad \mbox{\bf 1}_B\cirk f=f\cirk \mbox{\bf 1}_A=f, \\
(cat\: 2)\quad h\cirk (g\cirk f)=(h\cirk g)\cirk f,
\end{array}
\]
\[
\begin{array}{l}
(\times 1)\; \mbox{\bf 1}_A\times\mbox{\bf 1}_B=\mbox{\bf 1}_{A\times B}, \\%
[.05cm] (\times 2)\; (g_1\cirk g_2)\times(f_1\cirk f_2)=
(g_1\times
f_1)\cirk(g_2\times f_2), \\[.05cm]
(k^i)\; k^i_{B_1,B_2}\cirk(f_1\times f_2)=f_i\cirk k^i_{A_1,A_2}, \\[.05cm]
(w)\; w_B\cirk f=(f\times f)\cirk w_A, \\[.05cm]
(kw1)\; k^i_{A,A}\cirk w_A=\mbox{\bf 1}_A, \\[.05cm]
(kw2)\; (k^1_{A,B}\times k^2_{A,B})\cirk w_{A\times B}= \mbox{\bf 1}%
_{A\times B}, \\[.05cm]
(k)\quad {\mbox{\rm for }} f:A\rightarrow \mbox{\rm I},\; f=k_A, \\[.05cm]
(k\mbox{\rm O})\; k^1_{\mbox{\scriptsize{\rm O}},\mbox{\scriptsize{\rm O}}%
}=k^2_{\mbox{\scriptsize{\rm O}},\mbox{\scriptsize{\rm O}}},
\end{array}
\]
\[
\begin{array}{l}
(+1)\; \mbox{\bf 1}_A+\mbox{\bf 1}_B=\mbox{\bf 1}_{A+B}, \\[.05cm]
(+2)\; (g_1\cirk g_2)+(f_1\cirk f_2)=(g_1+f_1)\cirk(g_2+f_2), \\[.05cm]
(l^i)\; (f_1+f_2)\cirk l^i_{A_1,A_2}=l^i_{B_1,B_2}\cirk f_i, \\[.05cm]
(m)\; f\cirk m_A=m_B\cirk (f+f), \\[.05cm]
(lm1)\; m_A\cirk l^i_{A,A}=\mbox{\bf 1}_A, \\[.05cm]
(lm2)\; m_{A+B}\cirk (l^1_{A,B}+l^2_{A,B})=\mbox{\bf 1}_{A+B}, \\[.05cm]
(l)\quad {\mbox{\rm for }} f:\mbox{\rm O}\rightarrow A,\; f=l_A, \\[.05cm]
(l\mbox{\rm I})\; l^1_{\mbox{\scriptsize{\rm I}},\mbox{\scriptsize{\rm I}}%
}=l^2_{\mbox{\scriptsize{\rm I}},\mbox{\scriptsize{\rm I}}}.
\end{array}
\]

If we omit the equations $(k\mbox{\rm O})$ and $(l\mbox{\rm I})$, we obtain
the free bicartesian category generated by ${\cal L}$.

The free {\em sesquicartesian category} ${\cal S}$\ generated by $%
{\cal L}$ has as objects the formulae of ${\cal P}_{\times,+, %
\mbox{\scriptsize{\rm O}}}$. In that case, terms in which $k_A$
occurs are absent, and the equations $(k)$ and $(l\mbox{\rm I})$
are missing. The remaining terms and equations are as in ${\cal
D}$.

We shall call terms of ${\cal D}$\ in which the letters $l$ and
$m$ don't occur {\em K-terms}. (That means there are no subterms
of {\em K}-terms of the form $l_A$, $l^i_{A,B}$ and $m_A$.) Terms
of ${\cal D}$\ in which the letters $k$ and $w$ don't occur will
be called {\em L}-terms.

\section{Cut Elimination}

For inductive proofs on the length of objects or terms, it is
useful to be able to name the arrows of a category by terms that
contain no composition. In this section we shall prove for ${\cal
D}$\ and ${\cal S}$\ a theorem analogous to cut-elimination
theorems of logic, which stem from Gentzen \cite {gen35}. This
theorem says that every term of these categories is equal to a
term of a special sort, which in an appropriate language would be
translated by a composition-free term. We shall call terms of this
special sort {\em cut-free Gentzen terms}.

We define first the following operations on terms of ${\cal D}$,
which we call {\em Gentzen operations}:
\[
\begin{array}{ll}
K^1_B f=_{def.} f\cirk k^1_{A,B}, & L^1_B f=_{def.} l^1_{A,B}\cirk f, \\%
[.05cm]
K^2_A f=_{def.} f\cirk k^2_{A,B}, & L^2_A f=_{def.} l^2_{A,B}\cirk f, \\%
[.05cm] \langle f, g\rangle =_{def.} (f\times g)\cirk w_C, &
[f,g]=_{def.}m_C\cirk (f+g).
\end{array}
\]

Starting from the identities and the terms $k_A$ and $l_A$, for every $A$ of
${\cal P}$, and closing under the Gentzen operations, we obtain the set of
{\em cut-free Gentzen terms} of ${\cal D}$. The {\em Gentzen terms of} $%
{\cal D}$\ are obtained by closing cut-free Gentzen terms under
composition.

It is easy to show that every term of ${\cal D}$\ is equal in
${\cal D}$\ to a Gentzen term, since we have the following
equations:
\[
\begin{array}{ll}
k^1_{A,B}= K^1_B\mbox{\bf 1}_A, & l^1_{A,B}=L^1_B\mbox{\bf 1}_A, \\[.05cm]
k^2_{A,B}= K^2_A\mbox{\bf 1}_B, & l^2_{A,B}=L^2_A\mbox{\bf 1}_B, \\[.05cm]
w_A=\langle\mbox{\bf 1}_A,\mbox{\bf 1}_A\rangle, & m_A=[\mbox{\bf 1}_A,%
\mbox{\bf 1}_A], \\[.05cm]
f\times g=\langle K^1_C f, K^2_A g\rangle, & f+g= [L^1_D f, L^2_B g].
\end{array}
\]

We need the following equations of ${\cal D}$:
\[
\begin{array}{ll}
(K1)\quad g\cirk K^i_A f=K^i_A(g\cirk f), & (L1)\quad L^i_A g
\cirk
f=L^i_A(g\cirk f), \\[.05cm]
(K2)\quad K^i_Ag\cirk\langle f_1,f_2\rangle =g\cirk f_i, &
(L2)\quad [g_1,
g_2]\cirk L^i_A f=g_i\cirk f, \\[.05cm]
(K3)\quad \langle g_1,g_2\rangle \cirk f=\langle g_1\cirk f,
g_2\cirk f\rangle, & (L3)\quad g\cirk[f_1,f_2]=[g\cirk f_1, g\cirk
f_2],
\end{array}
\]
\noindent in order to prove the following theorem for ${\cal D}$.\\[.3cm]
C{\footnotesize UT} E{\footnotesize LIMINATION}.\quad {\em Every term is
equal to a cut-free Gentzen term.}

\vspace{2ex}

\noindent {\it Proof.} We first find for an arbitrary term of
${\cal D}$\ a Gentzen term $h$ equal to it. Let the {\em degree}
of a Gentzen term be the number of occurrences of Gentzen
operations in this term. Take a subterm $g\cirk f$ of $h$ such
that both $f$ and $g$ are cut-free Gentzen terms. We call such a term a {\em %
topmost cut}. We show that $g\cirk f$ is either equal to a
cut-free Gentzen term, or it is equal to a Gentzen term whose
topmost cuts are of strictly smaller degree than the degree of
$g\cirk f$. The possibility of eliminating the main compositions
of topmost cuts, and hence of finding for $h$ a cut-free Gentzen
term, follows by induction on degree.

The cases where $f$ or $g$ is $\mbox{\bf 1}_A$, or $f$ is $l_A$, or $g$ is $%
k_A$, are taken care of by $(cat\: 1)$, $(l)$ and $(k)$. The cases where $f$
is $K^i_A f^{\prime}$ or $g$ is $L^i_A g^{\prime}$ are taken care of by $%
(K1) $ and $(L1)$. And the cases where $f$ is $[f_1, f_2]$ or $g$ is $%
\langle g_1,g_2\rangle$ are taken care of by $(L3)$ and $(K3)$.

The following cases remain. If $f$ is $k_A$, then $g$ is of a form covered
by cases we dealt with above.

If $f$ is $\langle f_1,f_2\rangle$, then $g$ is either of a form covered by
cases above, or $g$ is $K^i_A g^{\prime}$, in which case we apply $(K2)$.

If $f$ is $L^i_A f^{\prime}$, then $g$ is either of a form covered
by cases above, or $g$ is $[g_1,g_2]$, in which case we apply
$(L2)$. This covers all possible cases. \qed

\vspace{2ex}

This proof, with the cases involving $k_A$ omitted, suffices to demonstrate
Cut Elimination for ${\cal S}$.

Let the {\em cut-free Gentzen K-terms} of ${\cal D}$\ be obtained
from the identities and the terms $k_A$ by closing under $+$ and
the Gentzen
operations $K^i$ and $\langle, \rangle$. The {\em Gentzen K-terms} of ${\cal %
D}$\ are obtained by closing the cut-free Gentzen {\em K}-terms
under composition. Let, dually, the {\em cut-free Gentzen L-terms}
of ${\cal D}$\ be obtained from the identities and the terms $l_A$
by closing under $\times$
and the Gentzen operations $L^i$ and $[,]$. The {\em Gentzen L-terms} of $%
{\cal D}$\ are obtained by closing the cut-free Gentzen {\em
L}-terms under composition.

Then we can prove the following version of Cut Elimination for the {\em K}%
-terms and the {\em L}-terms of ${\cal D}$.\\[.3cm]
C{\footnotesize UT} E{\footnotesize LIMINATION FOR} {\em K}-T{\footnotesize %
ERMS AND} {\em L}-T{\footnotesize ERMS}. \quad {\em Every K-term is equal to
a cut-free Gentzen K-term, and every L-term is equal to a cut-free Gentzen
L-term.}

\vspace{2ex}

\noindent {\it Proof.} It is easy to see that every {\em K}-term
is equal in ${\cal D}$\ to a Gentzen {\em K}-term. That this
Gentzen {\em K}-term is equal to a cut-free Gentzen {\em K}-term
is demonstrated as in the proof of Cut Elimination above by
induction on the degree of topmost cuts. We have to consider the
following additional cases.

If $f$ is $k_A$ or $\langle f_1,f_2\rangle$, then $g$ cannot be of the form $%
g_1+g_2$. If $f$ is $f_1+f_2$, and $g$ is not of a form already covered by
cases in the proof above, then $g$ is of the form $g_1+g_2$, in which case
we apply $(+2)$. This covers all possible cases.

Cut Elimination for {\em L}-terms follows by duality. \qed

\vspace{2ex}

Let the {\em cut-free Gentzen K-terms} of ${\cal S}$\ be terms of ${\cal S}$%
\ obtained from the identities by closing under $+$ and the Gentzen
operations $K^i$ and $\langle, \rangle$. The {\em Gentzen K-terms} of ${\cal %
S}$\ are obtained by closing the cut-free Gentzen {\em K}-terms of ${\cal S}$%
\ under composition. Let the {\em cut-free Gentzen L-terms} of ${\cal S}$\
be terms of ${\cal S}$\ obtained from the identities and the terms $l_A$ by
closing under $\times$ and the Gentzen operations $L^i$ and $[,]$. The {\em %
Gentzen L-terms} of ${\cal S}$\ are obtained by closing the cut-free Gentzen
{\em L}-terms of ${\cal S}$\ under composition.

Then we can establish Cut Elimination for {\em K}-terms and {\em L}-terms of
${\cal S}$, where equality in ${\cal D}$\ is replaced by equality in ${\cal S%
}$. We just proceed as in the proof above, with inapplicable cases involving
$k_A$ omitted.

\section{The Graphical Category}

We shall now define a graphical category ${\cal G}$ into which
${\cal D}$\ and ${\cal S}$\ can be mapped. The objects of ${\cal
G}$ are finite ordinals. An arrow $f:n\rightarrow m$ of ${\cal G}$
will be a binary relation from $n$ to $m$, i.e. a subset of
$n\times m$ with domain $n$ and codomain $m$. The identity
$\mbox{\bf 1}_n:n\rightarrow n$ of ${\cal G}$ is the identity
relation on $n$, and composition of arrows is composition of
relations.

For an object $A$ of ${\cal D}$, let $|A|$ be the number of
occurrences of
propositional letters in $A$. For example, $|(p\times(q+p))+(\mbox{\rm I}%
\times p)|$ is 4.

We now define a functor $G$ from ${\cal D}$\ to ${\cal G}$ such that $%
G(A)=|A|$. It is clear that $G(A\times B)=G(A+B)=|A|+|B|$. We define $G$ on
arrows inductively:
\[
\begin{array}{lcl}
G(\mbox{\bf 1}_A) & = & \{(x,x):x\in |A|\} = \mbox{\bf 1}_{|A|}, \\[.3cm]
G(k^1_{A,B}) & = & \{(x,x):x\in |A|\}, \\[.05cm]
G(k^2_{A,B}) & = & \{(x+|A|,x):x\in |B|\}, \\[.05cm]
G(w_A) & = & \{(x,x):x\in|A|\} \cup \{(x,x+|A|):x\in|A|\}, \\[.05cm]
G(k_A) & = & \emptyset,
\end{array}
\]
\[
\begin{array}{lcl}
G(l^1_{A,B}) & = & \{(x,x):x\in |A|\}, \\[.05cm]
G(l^2_{A,B}) & = & \{(x,x+|A|):x\in |B|\}, \\[.05cm]
G(m_A) & = & \{(x,x) :x\in|A|\} \cup \{(x+|A|,x):x\in |A|\}, \\[.05cm]
G(l_A) & = & \emptyset, \\[.3cm]
G(g\cirk f) & = & G(g) \cirk G(f),
\end{array}
\]
and for $f:A\rightarrow B$ and $g:C\rightarrow D$,
\[
G(f\times g)=G(f+g)=G(f) \cup \{(x+|A|,y+|B|):(x,y)\in G(g)\}.
\]

Though $G(\mbox{\bf 1}_A)$, $G(k^1_{A,B})$ and $G(l^1_{A,B})$ are the same
as sets of ordered pairs, in general they have different domains and
codomains, the first being a subset of $|A|\times |A|$, the second a subset
of $(|A|+|B|)\times |A|$, and the third a subset of $|A|\times(|A|+|B|)$. We
have an analogous situation in some other cases.

The arrows $G(f)$ of ${\cal G}$ can easily be represented graphically, by
drawings linking propositional letters, as it is illustrated in \cite{dos00b}%
. This is why we call this category "graphical".

It is easy to check that $G$ is a functor from ${\cal D}$\ to
${\cal G}$. We show by induction on the length of derivation that
if $f=g$ in ${\cal D}$, then $G(f)=G(g)$ in ${\cal G}$. (Of
course, $G$ preserves identities and composition.)

For the bicartesian structure of ${\cal G}$ we have that the operations $%
\times$ and $+$ on objects are both addition of ordinals, the operations $%
\times$ and $+$ on arrows coincide and are defined by the clauses for $%
G(f\times g)$ and $G(f+g)$, and the terminal and the initial
object also coincide: they are both the ordinal zero. The category
${\cal G}$ has zero arrows---namely, the empty relation. The
bicartesian category ${\cal G}$ is a linear category in the sense
of \cite{law97} (see p. 279). The functor $G$ from ${\cal D}$\ to
${\cal G}$ is not just a functor, but a bicartesian
functor; namely, a functor that preserves the bicartesian structure of $%
{\cal D}$.

We also have a functor defined analogously to $G$, which we call $G$ too,
from ${\cal S}$\ to ${\cal G}$. It is obtained from the definition of $G$
above by just rejecting clauses that are no longer applicable.

Our aim is to show that the functor $G$ from ${\cal S}$\ is
faithful.

\section{\emph{K-L} Normalization}

We shall say for a term of ${\cal D}$\ of the form $f_n\cirk \ldots\cirk f_1$%
, for some $n\geq 1$, where $f_i$ is composition-free, that it is {\em %
factorized}. By using $(\times 2)$, $(+2)$ and $(cat\: 1)$ it is easy to
show that every term of ${\cal D}$\ is equal to a factorized term of ${\cal D%
}$. A subterm $f_i$ in a factorized term $f_n\cirk \ldots\cirk
f_1$ is called a {\em factor}.

A term of ${\cal D}$\ where all the atomic terms are identities
will be
called a {\em complex identity}. According to $(\times 1)$, $(+1)$ and $%
(cat\: 1)$, every complex identity is equal to an identity. A factor which
is a complex identity will be called an {\em identity factor}. It is clear
that if $n>1$, we can omit in a factorized term every identity factor, and
obtain a factorized term equal to the original one.

A term of ${\cal D}$\ is said to be in {\em K-L normal form} iff
it is of
the form $g\cirk f:A\rightarrow B$ for $f$ a {\em K}-term and $g$ an {\em L}%
-term. Note that {\em K-L} normal forms are not unique, since
$(m_A\times m_A)\cirk w_{A+A}$ and $m_{A\times A}\cirk(w_A+w_A)$,
which are both equal to $w_A\cirk m_A$, are both in {\em K-L}
normal form.

We can prove the following proposition for ${\cal D}$.\\[.3cm]
{\em K-L} N{\footnotesize ORMALIZATION}. \quad {\em Every term is equal to a
term in} {\em K-L} {\em normal form}.

\vspace{2ex}

\noindent {\it Proof.} Suppose $f:B\rightarrow C$
is a composition-free {\em K}-term that is not a complex identity, and $%
g:A\rightarrow B$ is a composition-free {\em L}-term that is not a
complex identity. We show by induction on the length of $f\cirk g$
that
\[
(*)\quad f\cirk g=g^{\prime}\cirk f^{\prime}\; {\mbox{\rm or }}
f\cirk g=f^{\prime}\; {\mbox{\rm or }} f\cirk g=g^{\prime}
\]
for $f^{\prime}$ a composition-free {\em K}-term and $g^{\prime}$ a
composition-free {\em L}-term.

We shall not consider below cases where $g$ is $m_B$, which are easily taken
care of by $(m)$. Cases where $f$ is $k_B$ or $g$ is $l_B$ are easily taken
care of by $(k)$ and $(l)$. The following cases remain.

If $f$ is $k^i_{C,E}$ and $g$ is $g_1\times g_2$, then we use $(k^i)$. If $f$
is $w_B$, then we use $(w)$. If $f$ is $f_1\times f_2$ and $g$ is $g_1\times
g_2$, then we use $(\times 2)$, the induction hypothesis, and perhaps $%
(cat\: 1)$.

Finally, if $f$ is $f_1+f_2$, then we have the following cases. If $g$ is $%
l^i_{B_1,B_2}$, then we use $(l^i)$. If $g$ is $g_1+g_2$, then we use $(+2)$%
, the induction hypothesis, and perhaps $(cat\: 1)$. This proves $(*)$.

Every term of ${\cal D}$\ is equal to an identity or to a factorized term $%
f_n\cirk\ldots\cirk f_1$ without identity factors. Every factor $f_i$ of $%
f_n\cirk\ldots\cirk f_1$ is either a {\em K}-term or an {\em L}-term or, by $%
(cat\: 1)$, $(\times2)$ and $(+2)$, it is equal to
$f_i^{\prime\prime}\cirk f_i^{\prime}$ where $f_i^{\prime}$ is a
{\em K}-term and $f_i^{\prime\prime}$ is an {\em L}-term. For
example, $(k^1_{A,B}\times l^1_{C,D})+(w_E + l_F)$ is equal to
\[
((\mbox{\bf 1}_A\times l^1_{C,D})+(\mbox{\bf 1}_{E\times E}) +
l_F)\cirk
((k^1_{A,B}\times\mbox{\bf 1}_C)+(w_E + \mbox{\bf 1}_{%
\mbox{\scriptsize{\rm
O}}})).
\]
Then it is clear that by applying $(*)$ repeatedly, and by applying perhaps $%
(cat\: 1)$, we obtain a term in {\em K-L} normal form. \qed

\vspace{2ex}

Note that to reduce a term of ${\cal D}$\ to {\em K-L} normal form
we have used in this proof all the equations of ${\cal D}$\ except
$(kw1)$, $(kw2)$, $(lm1)$, $(lm2)$, $(k\mbox{\rm O})$ and
$(l\mbox{\rm I})$.

The definition of {\em K-L} normal form for ${\cal S}$\ is the
same. Then the proof above, with some parts omitted, establishes
{\em K-L} Normalization also for ${\cal S}$.

\section{Coherence for Sesquicartesian Categories}

We shall prove in this section that the functor $G$ from ${\cal S}$\ to $%
{\cal G}$ is faithful, i.e. we shall show that we have coherence
for sesquicartesian categories. These categories are interesting
because the category {\bf Set} of sets with functions, where
cartesian product is $\times$, disjoint union is $+$, and the
empty set is O, is such a category. As a matter of fact, every
bicartesian closed category is a sesquicartesian category. A
bicartesian closed category is a bicartesian category that is
cartesian closed, i.e., every functor $A\times ...$ has a right
adjoint $...^{A}$. And in every cartesian closed category with an
initial object O we have $(k\mbox{\rm O})$, because $Hom(\mbox{\rm
O} \times \mbox{\rm O}, \mbox{\rm O})\cong Hom(\mbox{\rm O},
\mbox{\rm O}^{\mbox{\scriptsize{\rm O}}})$. In every cartesian
category in which $(k\mbox{\rm O})$ holds we have that
$Hom(A,\mbox{\rm O})$ is a singleton or empty, because for
$f,g:A\rightarrow \mbox{\rm O}$ we have
$k^{1}_{\mbox{\scriptsize{\rm O}},\mbox{\scriptsize{\rm O}}} \cirk
\langle f,g \rangle = k^{2}_{\mbox{\scriptsize{\rm
O}},\mbox{\scriptsize{\rm O}}} \cirk \langle f,g \rangle$ (cf.
\cite{LS86}, Proposition 8.3, p. 67).

The category {\bf Set} shows that sesquicartesian categories are
not maximal in the following sense. In the category ${\cal S}$\
the
equations $l_{\mbox{\scriptsize{\rm O}}\times A}\cirk k_{\mbox{\scriptsize{%
\rm O}},A}^{1}=\mbox{\bf 1}_{\mbox{\scriptsize{\rm O}}\times A}$ and $%
l_{A\times \mbox{\scriptsize{\rm O}}}\cirk k_{A,\mbox{\scriptsize{\rm O}}%
}^{2}=\mbox{\bf 1}_{A\times \mbox{\scriptsize{\rm O}}}$ (in which only terms
of ${\cal S}$\ occur) don't hold, but they hold in {\bf Set}, and {\bf Set}
is not a preorder. That these two equations don't hold in ${\cal S}$\
follows from the fact that $G$ is a functor from ${\cal S}$\ to ${\cal G}$,
but $G(l_{\mbox{\scriptsize{\rm O}}\times p}\cirk k_{%
\mbox{\scriptsize{\rm
O}},p}^{1})$ and $G(l_{p\times \mbox{\scriptsize{\rm O}}}\cirk k_{p,%
\mbox{\scriptsize{\rm O}}}^{2})$ are empty, whereas $G(\mbox{\bf 1}_{%
\mbox{\scriptsize{\rm O}}\times p})$ and $G(\mbox{\bf 1}_{p\times %
\mbox{\scriptsize{\rm O}}})$ contain $(0,0)$. In the case of
cartesian categories and categories with binary products and sums,
we had coherence, and used this coherence to prove maximality in
\cite{dos00} and \cite{dos00b}. With sesquicartesian categories,
however, coherence and maximality don't go hand in hand any more.

First we prove the following lemmata. \\[.3cm]
{L{\footnotesize EMMA} 6.1.}\hspace{1em} {\em A constant object of} ${\cal P}%
_{\times,+,\mbox{\scriptsize{\rm O}}}$ {\em is isomorphic} {\em in} ${\cal S}
$\ {\em to} O.

\vspace{2ex}

\noindent {\it Proof.} We have in ${\cal S}$\ the isomorphisms
\[
\begin{array}{ll}
k^1_{\mbox{\scriptsize{\rm O}},\mbox{\scriptsize{\rm O}}}=k^2_{%
\mbox{\scriptsize{\rm O}},\mbox{\scriptsize{\rm O}}}:\mbox{\rm O}\times%
\mbox{\rm O}\rightarrow\mbox{\rm O},{\mbox{\hspace{3em}}} & l_{%
\mbox{\scriptsize{\rm O$\times$O}}}=w_{\mbox{\scriptsize{\rm O}}}:\mbox{\rm O}%
\rightarrow\mbox{\rm O}\times\mbox{\rm O}, \\[.05cm]
m_A\cirk(\mbox{\bf 1}_A +l_A):A+\mbox{\rm O}\rightarrow A, & l^1_{A,%
\mbox{\scriptsize{\rm O}}}:A\rightarrow A+\mbox{\rm O}, \\[.05cm]
m_A\cirk(l_A +\mbox{\bf 1}_A):\mbox{\rm O} +A\rightarrow A, & l^2_{%
\mbox{\scriptsize{\rm O}},A}:A\rightarrow \mbox{\rm O}
+A.\quad\dashv
\end{array}
\]

\vspace{2ex}

\noindent {L{\footnotesize EMMA} 6.2.}\hspace{1em} {\em If} $f,g:A\rightarrow B$ {\em %
are terms of} ${\cal S}$\ {\em and either} $A$ {\em or} $B$ {\em is
isomorphic in} ${\cal S}$\ {\em to} O, {\em then} $f=g$ {\em in} ${\cal S}$.

\vspace{2ex}

\noindent {\it Proof.}
Suppose $i: A\rightarrow %
\mbox{\rm O}$ is an isomorphism in ${\cal S}$. Then from
\[
f\cirk i^{-1}=g\cirk i^{-1}=l_B
\]
we obtain $f=g$.

Suppose $i:B\rightarrow \mbox{\rm O}$ is an isomorphism in ${\cal S}$. Then
from
\[
k^1_{\mbox{\scriptsize{\rm O}},\mbox{\scriptsize{\rm
O}}}\cirk\langle i\cirk
f,i\cirk g\rangle= k^2_{\mbox{\scriptsize{\rm O}},\mbox{\scriptsize{\rm O}}%
}\cirk\langle i\cirk f,i\cirk g\rangle
\]
we obtain $i\cirk f=i\cirk g$, which yields $f=g$. \qed

\vspace{2ex}

\noindent {L{\footnotesize EMMA} 6.3.}\hspace{1em} {\em If} $f,g:A\rightarrow B$ {\em %
are terms of} ${\cal S}$\ {\em and} $G(f)=G(g)=\emptyset$ {\em in} ${\cal G}$%
, {\em then} $f=g$ {\em in} ${\cal S}$.

\vspace{2ex}

\noindent {\it Proof.} By {\em K}-$L$ Normalization, $f=f_2\cirk
f_1$ for $f_1:A\rightarrow C$ a
{\em K}-term and $f_2:C\rightarrow B$ an {\em L}-term. Since for every $%
z\in|C|$ there is a $y\in|B|$ such that $(z,y)\in G(f_2)$, we must have $%
G(f_1)$ empty; otherwise, $G(f)$ would not be empty. On the other hand, if
there is a propositional letter in $C$ at the $z$-th place, there must be
for some $x\in |A|$ a pair $(x,z)$ in $G(f_1)$. So $C$ is a constant object
of ${\cal P}_{\times,+,\mbox{\scriptsize{\rm O}}}$.

Analogously, $g=g_2\cirk g_1$ for $g_1:A\rightarrow D$ a {\em K}-term, and $%
g_2:D\rightarrow B$ an {\em L}-term. As before, $D$ is a constant object of $%
{\cal P}_{\times,+,\mbox{\scriptsize{\rm O}}}$. By Lemma 6.1,
there is an isomorphism $i:C\rightarrow D$ of ${\cal S}$, and
$f=f_2\cirk i^{-1}\cirk
i\cirk f_1$. By Lemmata 6.1 and 6.2, we obtain $i\cirk f_1=g_1$ and $%
f_2\cirk i^{-1}=g_2$, from which $f=g$ follows. \qed

\vspace{2ex}

We shall next prove the following coherence proposition.\\[.3cm]
S{\footnotesize ESQUICARTESIAN} C{\footnotesize OHERENCE}. \quad
{\em If} $f,g:A\rightarrow B$ {\em are terms of} ${\cal S}$\ {\em and} $%
G(f)=G(g)$ {\em in} ${\cal G}$, {\em then} $f=g$ {\em in} ${\cal S}$.

\vspace{2ex}

\noindent {\it Proof.} Lemma 6.3 covers the case when
$G(f)=G(g)=\emptyset$. So we assume $G(f)=G(g)\neq\emptyset$, and
proceed by induction on the sum of the lengths of $A$ and $B$. In
this induction we need not consider the cases when either $A$ or
$B$ is a constant object; otherwise, $G(f)$ and $G(g)$ would be
empty. So in the basis of the induction, when both of $A$ and $B$
are atomic, we consider only the case when both of $A$ and $B$ are
propositional letters. In this case we conclude by Cut Elimination
that $f$ and $g$ exist iff $A$ and $B$ are the same propositional
letter $p$, and $f=g=\mbox{\bf 1}_p$ in ${\cal S}$. (We could
conclude the same thing by interpreting ${\cal S}$\ in
conjunctive-disjunctive logic.) Note that we didn't need here the
assumption $G(f)=G(g)$.

If $A$ is $A_1+A_2$, then $f\cirk l^1_{A_1,A_2}$ and $g\cirk
l^1_{A_1,A_2}$ are of type $A_1\rightarrow B$, while $f\cirk
l^2_{A_1,A_2}$ and $g\cirk l^2_{A_1,A_2}$ are of type
$A_2\rightarrow B$. We also have
\[
\begin{array}{rcl}
G(f\cirk l^i_{A_1,A_2}) & = & G(f)\cirk G(l^i_{A_1,A_2}) \\
& = & G(g)\cirk G(l^i_{A_1,A_2}) \\
& = & G(g\cirk l^i_{A_1,A_2}),
\end{array}
\]
whence, by the induction hypothesis, or Lemma 6.3,
\[
f\cirk l^i_{A_1,A_2}=g\cirk l^i_{A_1,A_2}
\]
in ${\cal S}$. Then we infer that
\[
[f\cirk l^1_{A_1,A_2},f\cirk l^2_{A_1,A_2}]=[g\cirk
l^1_{A_1,A_2},g\cirk l^2_{A_1,A_2}],
\]
from which it follows that $f=g$ in ${\cal S}$. We proceed analogously if $B$
is $B_1\times B_2$.

Suppose now $A$ is $A_1\times A_2$ or a propositional letter, and $B$ is $%
B_1+B_2$ or a propositional letter, but $A$ and $B$ are not both
propositional letters. Then, by Cut Elimination, $f$ is equal in
${\cal S}$ either to a term of the form $f^{\prime}\cirk
k^i_{A_1,A_2}$, or to a term of the form $l^i_{B_1,B_2}\cirk
f^{\prime}$. Suppose $f=f^{\prime}\cirk k^1_{A_1,A_2}$. Then for
every $(x,y)\in G(f)$ we have $x\in |A_1|$. (We reason analogously
when $f=f^{\prime}\cirk k^2_{A_1,A_2}$.)

By Cut Elimination too, $g$ is equal in ${\cal S}$ either to a
term of the form ${g'\cirk k^i_{A_1,A_2}}$, or to a term of the
form ${l^i_{B_1,B_2}\cirk g'}$. In the first case we must have
${g=g'\cirk k^1_{A_1,A_2}}$, because ${G(g)=G(f^{\prime}\cirk
k^1_{A_1,A_2})\neq\emptyset}$, and then we apply the induction
hypothesis to derive ${f'=g'}$ from ${G(f')=G(g')}$. Hence ${f=g}$
in $\cal S$.

Suppose ${g=l^1_{B_1,B_2}\cirk g'}$. (We reason analogously when
${g=l^2_{B_1,B_2}\cirk g'}$.) Let ${f'':A_1\rightarrow B_1+B_2''}$
be the substitution instance of ${f':A_1\rightarrow B_1+B_2}$
obtained by replacing every occurrence of propositional letter in
$B_2$ by O. There is an isomorphism ${i:B_2''\rightarrow\mbox{\rm
O}}$ by Lemma 6.1, and $f''$ exists because in $G(f)$, which is
equal to ${G(l^1_{B_1,B_2}\cirk g')}$, there is no pair ${(x,y)}$
with ${y\geq|B_1|}$. So we have an arrow ${f''':A_1\rightarrow
B_1}$, which we define as $[\mbox{\bf
1}_{B_1},l_{B_1}]\cirk(\mbox{\bf 1}_{B_1}+ i)\cirk f''$. It is
easy to verify that ${G(l^1_{B_1,B_2}\cirk f''')=G(f')}$, and that
${G(f'''\cirk k^1_{A_1,A_2})=G(g')}$. By the induction hypothesis,
we obtain ${l^1_{B_1,B_2}\cirk f'''=f'}$ and ${f'''\cirk
k^1_{A_1,A_2}=g'}$, from which we derive ${f=g}$.

We reason analogously when ${f=l^i_{B_1,B_2}\cirk f'}$. \qed

\vspace{2ex}

To verify whether for $f,g:A\rightarrow B$ in the language of ${\cal S}$\ we
have $f=g$ in ${\cal S}$\ it is enough to draw $G(f)$ and $G(g)$, and check
whether they are equal, which is clearly a finite task. So we have here an
easy decision procedure for the equations of ${\cal S}$.

It is clear that we also have coherence for coherent dual sesquicartesian
categories, which are categories with arbitrary finite products and nonempty
finite sums where $(l\mbox{\rm I})$ holds.

As a consequence of Cut Elimination, of the functoriality of $G$ from ${\cal %
D}$, and of Sesquicartesian Coherence, we obtain that ${\cal S} $\
is a full subcategory of ${\cal D}$.

\section{Restricted Coherence for Dicartesian Catego\-ries}

First we give a few examples of dicartesian categories. Note that
the bicartesian category {\bf Set}, where product, sum and the
initial object are taken as usual (see the beginning of the
previous section), and a singleton is the
terminal object I, is not a dicartesian category. The equation $(l%
\mbox{\rm I})$ does not hold in {\bf Set}.

Every bicartesian category in which the terminal and the initial
object are isomorphic is a dicartesian category. Such is, for
instance, the category {\bf Set}* of pointed sets, i.e. sets with
a distinguished element $*$ and $*$-preserving maps, which is
isomorphic to the category of sets with partial maps. In {\bf
Set}* the objects I and O are both $\{*\}$, the product $\times $
is cartesian product, the sum $A+B$ of the objects $A$ and $B$ is
$\{(a,*):a\in A\}\cup \{(*,b):b\in B\}$, with $(*,*)$ being the
$*$ of the product and of the sum, while
\begin{eqnarray*}
&&l_{A,B}^{1}(a)=(a,*),\quad l_{B,A}^{2}(a)=(*,a), \\
&&m_{A}(a,b)=\left\{
\begin{array}{ll}
a & \mbox{\rm{if} }a\neq * \\
b & \mbox{\rm{if} }b\neq * \\
\ast & \mbox{\rm{if} }a=b=*
\end{array}
\right. \\
&&\,(f+g)(a,b)=(f(a),g(b)).
\end{eqnarray*}
In {\bf Set}* we have that I$+$I is isomorphic to I.

A fortiori, every bicartesian category in which all finite
products and sums are isomorphic (i.e. every linear category in
the sense of \cite{law97}, p. 279), is a dicartesian category.
Such are, for example, the category of commutative monoids with
monoid homomorphisms, and its subcategory of vector spaces over a
fixed field with linear transformations. We shall next present a
concrete dicartesian category in which the terminal and the
initial object, as well as finite products and sums in general,
are not isomorphic.

As sesquicartesian categories were not maximal, so dicartesian
categories are not maximal either. This is shown by the category
{\bf Set}*$+\emptyset $, which is the category of pointed sets
extended with the empty set $\emptyset $
as an additional object. The arrows of {\bf Set}*$+\emptyset $ are the $*$%
-preserving maps plus the empty maps with domain $\emptyset $. Let I be $%
\{*\}$, let O be $\emptyset $, let $\times $ be cartesian product,
and let the sum $A+B$ be defined as in {\bf Set}*, with the same
clause. (If both $A$ and $B$ are $\emptyset $, then $A+B$ is
$\emptyset $.) If $A \times B$ and $A+B$ are not $\emptyset $,
then their $*$ is $(*,*)$, as in {\bf Set}*. Next, let $l_{B}$ be
$\emptyset :\emptyset \rightarrow B$, and if $A$ is not $\emptyset
$, let $l_{A,B}^{1}$, $l_{B,A}^{2}$ and $m_{A}$ be defined as in
{\bf Set}*. If $A$ is $\emptyset $, then $l_{A,B}^{1}$ is
$\emptyset :\emptyset \rightarrow \emptyset +B$, $l_{B,A}^{2}$ is
$\emptyset :\emptyset \rightarrow B+\emptyset $ and $m_{A}$ is
$\emptyset :\emptyset \rightarrow \emptyset $. We also have
\begin{eqnarray*}
&&\,(f+g)(a,b)=\left\{
\begin{array}{ll}
(f(a),g(b)) & {\mbox{\rm{if} neither }}f\;{\mbox{\rm nor }}g\;
{\mbox{\rm is an empty map}} \\
(f(a),*) & {\mbox{\rm {if} }}f\;{\mbox{\rm is not an empty map and }}g\;
{\mbox{\rm is}} \\
(*,g(b)) & {\mbox{\rm {if} }}g\;{\mbox{\rm is not an empty map and }}f\;
{\mbox{\rm is.}}
\end{array}
\right.
\end{eqnarray*}
If $f$ and $g$ are both empty maps, then $f+g$ is the empty map
with appropriate codomain.

With other arrows and operations on arrows defined in the obvious
way, we can check that {\bf Set}*$+\emptyset $ is a dicartesian
category. In this category we have that $\emptyset \times
A=A\times \emptyset =\emptyset $. In the category ${\cal D}$ the
equations $l_{\mbox{\scriptsize{\rm O}}\times A}\cirk k_{\mbox{\scriptsize{%
\rm O}},A}^{1}=\mbox{\bf 1}_{\mbox{\scriptsize{\rm O}}\times A}$ and $%
l_{A\times \mbox{\scriptsize{\rm O}}}\cirk k_{A,\mbox{\scriptsize{\rm O}}%
}^{2}=\mbox{\bf 1}_{A\times \mbox{\scriptsize{\rm O}}}$ (in which
only terms of ${\cal D}$ occur) don't hold, as we explained at the
beginning of the
previous section, but they hold in {\bf Set}*$+\emptyset $, and {\bf Set}*$%
+\emptyset $ is not a preorder.

Note that a dicartesian category ${\cal C}$ is cartesian closed
only if $%
{\cal C}$ is a preorder. We have $Hom_{{\cal C}}(A,B)\cong Hom_{{\cal C}}($I$%
,B^{A})$, and for $f,g:~$I$~\rightarrow D$ with $(l$I$)$ we obtain $%
[f,g]\cirk l_{{\mbox{\scriptsize{\rm I}}},{\mbox{\scriptsize{\rm
I}}}}^{1}=[f,g]\cirk l_{{\mbox{\scriptsize{\rm
I}}},{\mbox{\scriptsize{\rm I}}}}^{2}$,
which gives $f=g$. The category {\bf Set} is cartesian closed, whereas {\bf %
Set}* and {\bf Set}*$+\emptyset $ are not.

We can prove the following statements, which extend Lemmata 6.1-6.3. \\[.3cm]
{L{\footnotesize EMMA} 7.1.}\hspace{1em} {\em A constant object
of} ${\cal P} $ {\em is isomorphic in} ${\cal D}$\ {\em to either}
O {\em or} I.

\vspace{2ex}

\noindent {\it Proof.} In addition to the isomorphisms of the
proof of Lemma 6.1, we have in ${\cal D}$\ the isomorphisms
\[
\begin{array}{ll}
k_{\mbox{\scriptsize{\rm I}}}=m_{\mbox{\scriptsize{\rm I}}}:\mbox{\rm I} +%
\mbox{\rm I}\rightarrow\mbox{\rm I},{\mbox{\hspace{3em}}} & l^1_{%
\mbox{\scriptsize{\rm I}},\mbox{\scriptsize{\rm I}}}=l^2_{%
\mbox{\scriptsize{\rm I}},\mbox{\scriptsize{\rm I}}}:\mbox{\rm I}\rightarrow%
\mbox{\rm I} +\mbox{\rm I}, \\[.05cm]
k^1_{A,\mbox{\scriptsize{\rm I}}}:A\times \mbox{\rm I}\rightarrow A, & (%
\mbox{\bf 1}_A\times k_A)\cirk w_A:A\rightarrow A\times \mbox{\rm I}, \\%
[.05cm] k^2_{\mbox{\scriptsize{\rm I}},A}:\mbox{\rm I}\times
A\rightarrow A, & (k_A\times\mbox{\bf 1}_A)\cirk w_A:A\rightarrow
\mbox{\rm I}\times A.\quad \dashv
\end{array}
\]

\vspace{2ex}

\noindent {L{\footnotesize EMMA} 7.2.}\hspace{1em} {\em If} $f,g:A\rightarrow B$ {\em %
are terms of} ${\cal D}$\ {\em and either} $A$ {\em or} $B$ {\em
is
isomorphic in} ${\cal D}$\ {\em to} O {\em or} I, {\em then} $f=g$ {\em in} $%
{\cal D}$.

\vspace{2ex}

\noindent {\it Proof.} We repeat what we had in the proof of Lemma
6.2, and reason dually when $A$ or $B$ is isomorphic to I. \qed

\vspace{2ex}

\noindent {R{\footnotesize ESTRICTED} D{\footnotesize ICARTESIAN}
C{\footnotesize OHERENCE} I.}\hspace{1em} {\em If} $f,g:A\rightarrow B$ {\em %
are terms of} ${\cal D}$\ {\em and} $G(f)=G(g)=\emptyset$ {\em in}
${\cal G}$, {\em then} $f=g$ {\em in} ${\cal D}$.

\vspace{2ex}

\noindent {\it Proof.} As in the proof of Lemma
6.3, by {\em K}-$L$ Normalization, we have $f=f_2\cirk f_1$ for $%
f_1:A\rightarrow C$ a {\em K}-term, $f_2:C\rightarrow B$ an {\em L}-term,
and $C$ a constant object of ${\cal P}$; we also have $g=g_2\cirk g_1$ for $%
g_1:A\rightarrow D$ a {\em K}-term, $g_2:D\rightarrow B$ an {\em L}-term,
and $D$ a constant object of ${\cal P}$. Next we apply Lemma 7.1. If $C$ and
$D$ are both isomorphic to O, we reason as in the proof of Lemma 6.3, and we
reason analogously when they are both isomorphic to I. If $i:C\rightarrow %
\mbox{\rm O}$ and $j:\mbox{\rm I}\rightarrow D$ are isomorphisms
of ${\cal D} $, then we have
\[
\begin{array}{rcl}
f_2\cirk f_1 & = & g_2\cirk j\cirk k_{\mbox{\scriptsize{\rm
O}}}\cirk i\cirk
f_1,\; {\mbox{\rm by Lemma 7.2,}} \\[.05cm]
& = & g_2\cirk g_1,\; {\mbox{\rm by Lemma 7.2,}}
\end{array}
\]
and so $f=g$ in ${\cal D}$. (Note that $k_{\mbox{\scriptsize{\rm O}}}=l_{%
\mbox{\scriptsize{\rm I}}}$ in ${\cal D}$.) \qed

\vspace{2ex}

Besides Restricted Dicartesian Coherence I, we can prove another
partial coherence result for ${\cal D}$. For that result we need
the following lemma, and the definitions that follow.

\vspace{2ex}

\noindent {L{\footnotesize EMMA} 7.3\hspace{1em} {\em If ${u:
A\rightarrow B_1+B_2}$ is such that ${G(u)\neq\emptyset}$ and
there is no ${(x,y)}$ in $G(u)$ such that ${y\geq|B_1|}$ and $+$
does not occur in $A$, then there is an arrow ${v: A\rightarrow
B_1}$ such that ${G(u)=G(l^1_{B_1,B_2}\cirk v)}$}.

\vspace{2ex}

\noindent {\it Proof.} By induction on the length of $A$. Suppose
$u$ is a cut-free term. If $A$ is a propositional letter, then by
the assumption on ${G(u)}$ we have ${u=L^1_{B_2}u'}$ and we can
take ${v=u'}$.

If $A$  is not a propositional letter and $u$ is not of the form
${L^1_{B_2}u'}$ (by the assumption on ${G(u)}$, the term $u$
cannot be of the form ${L^2_{B_1}u'}$), then since $+$ does not
occur in $A$ we have that $u$ is of the form ${K^i_{A''}u'}$ for
${u':A'\rightarrow B_1+B_2}$. Since ${G(u')\neq\emptyset}$ and
there is no ${(x,y)}$ in $G(u')$ such that ${y\geq|B_1|}$ and $+$
does not occur in $A'$ we may apply the induction hypothesis to
$u'$ and obtain $v'$ such that ${G(u')=G(l^1_{B_1,B_2}\cirk v')}$
and hence we can take ${v=v'\cirk k^i}$. \qed

\vspace{2ex}

A formula $C$ of ${\cal P}$ is called a {\em contradiction} when
there is in ${\cal D}$ an arrow of the type ${C\rightarrow
\mbox{\rm O}}$. For every formula that is not a contradiction
there is a substitution instance isomorphic to I. Suppose $C$ is
not a contradiction, and let $C^{\mbox{\scriptsize{\rm I}}}$ be
obtained from $C$ by substituting I for every propositional
letter. If $C^{\mbox{\scriptsize{\rm I}}}$ were not isomorphic to
I, then by Lemma 7.1 we would have an isomorphism
${i:C^{\mbox{\scriptsize{\rm I}}}\rightarrow \mbox{\rm O}}$. Since
there is obviously an arrow ${u:C\rightarrow
C^{\mbox{\scriptsize{\rm I}}}}$ formed by using $k_p$, we would
have ${i\cirk u:C\rightarrow\mbox{\rm O}}$, and $C$ would be a
contradiction.

A formula $C$ of ${\cal P}$ is called a {\em tautology} when there
is in ${\cal D}$ an arrow of the type $\mbox{\rm I}\rightarrow C$.
For every formula that is not a tautology there is a substitution
instance isomorphic to O. (This is shown analogously to what we
had in the preceding paragraph.)

A formula of ${\cal P}$ is called O-{\em normal} when for every
subformula $D\times C$ or $C\times D$ of it with $C$ a
contradiction, there is no occurrence of $+$ in $D$. A formula of
${\cal P}$ is called I-{\em normal} when for every subformula
$D+C$ or $C+D$ of it with $C$ a tautology, there is no occurrence
of $\times$ in $D$.

We can now formulate our second partial coherence result for
dicartesian categories.

\vspace{2ex}

\noindent R{\footnotesize ESTRICTED} D{\footnotesize ICARTESIAN}
C{\footnotesize OHERENCE} II.\quad{\em If ${f,g:A\rightarrow B}$
are terms of ${\cal D}$ such that ${G(f)=G(g)}$ and either $A$ is
O-normal or $B$ is I-normal, then ${f=g}$ in ${\cal D}$}.

\vspace{2ex}

\noindent {\it Proof.} Suppose $A$ is O-normal. Restricted
Dicartesian Coherence I covers the case when
$G(f)=G(g)=\emptyset$. So we assume $G(f)=G(g)\neq\emptyset$, and
proceed as in the proof of Sesquicartesian Coherence by induction
on the sum of the lengths of $A$ and $B$. The basis of this
induction and the cases when $A$ is of the form ${A_1+A_2}$ or $B$
is of the form ${B_1\times B_2}$ are settled as in the proof of
Sesquicartesian Coherence.

Suppose $A$ is ${A_1\times A_2}$ or a propositional letter and $B$
is ${B_1+B_2}$ or a propositional letter, but $A$ and $B$ are not
both propositional letters. (The cases when $A$ or $B$ is a
constant object are excluded by the assumption that
$G(f)=G(g)\neq\emptyset$.) We proceed then as in the proof of
Sesquicartesian Coherence until we reach the case when ${f=f'\cirk
k^1_{A_1,A_2}}$ and ${g=l^1_{B_1,B_2}\cirk g'}$.

Suppose $A_2$ is not a contradiction. Then there is an instance
${A_2^{\mbox\scriptsize{\rm I}}}$ of $A_2$ and an isomorphism
${i:\mbox{\rm I}\rightarrow A_2^{\mbox\scriptsize{\rm I}}}$. (To
obtain ${A_2^{\mbox\scriptsize{\rm I}}}$ we substitute I for every
letter in $A_2$.) Let ${g'':A_1\times A_2^{\mbox\scriptsize{\rm
I}}\rightarrow B_1}$ be the substitution instance of
${g':A_1\times A_2\rightarrow B_1}$ obtained by replacing every
occurrence of propositional letter in $A_2$ by I. Such a term
exists because in ${G(g)}$, which is equal to ${G(f'\cirk
k^1_{A_1,A_2})}$, there is no pair ${(x,y)}$ with ${x\geq|A_1|}$.

So we have an arrow $g'''=g''\cirk(\mbox{\bf 1}_{A_1}\times
i)\cirk\langle\mbox{\bf 1}_{A_1},k_{A_1}\rangle:A_1\rightarrow
B_1$. It is easy to verify that $G(l^1_{B_1,B_2}\cirk g''')=G(f')$
and that $G(g'''\cirk k^1_{A_1,A_2})=G(g')$. By the induction
hypothesis we obtain $l^1_{B_1,B_2}\cirk g'''=f'$ and $g'''\cirk
k^1_{A_1,A_2}=g'$, from which we derive $f=g$.

Suppose $A_2$ is a contradiction. Then by the assumption that $A$
is O-normal we have that $+$ does not occur in $A_1$. We may apply
Lemma 7.3 to ${f':A_1\rightarrow B_1+B_2}$ to obtain
$f''':A_1\rightarrow B_1$ such that $G(f')=G(l^1_{B_1,B_2}\cirk
f''')$. It is easy to verify that then $G(g')=G(f'''\cirk
k^1_{A_1,A_2})$, and we may proceed as in the proof of
Sesquicartesian Coherence.

We proceed analogously when $B$ is I-normal, relying on a lemma
dual to Lemma 7.3.\qed

\vspace{2ex}

Let $A^0_{\mbox\scriptsize\rm O}$ be $A\times \mbox{\rm O}$, and
let $A^{n+1}_{\mbox\scriptsize\rm O}$ be
$(A^n_{\mbox\scriptsize\rm O}+\mbox{\rm I})\times \mbox{\rm O}$.
Let $f^0_{\mbox\scriptsize\rm O}$ be $f\times \mbox{\bf
1}_{\mbox\scriptsize\rm O}$, and let $f^{n+1}_{\mbox\scriptsize\rm
O}$ be $(f^n_{\mbox\scriptsize\rm O}+\mbox{\bf
1}_{\mbox\scriptsize\rm I})\times \mbox{\bf
1}_{\mbox\scriptsize\rm O}$. Let, dually,
$A^0_{\mbox\scriptsize\rm I}$ be $A+ \mbox{\rm I}$, and let
$A^{n+1}_{\mbox\scriptsize\rm I}$ be $(A^n_{\mbox\scriptsize\rm
I}\times\mbox{\rm O})+ \mbox{\rm I}$. Let
$f^0_{\mbox\scriptsize\rm I}$ be $f+ \mbox{\bf
1}_{\mbox\scriptsize\rm I}$, and let $f^{n+1}_{\mbox\scriptsize\rm
I}$ be $(f^n_{\mbox\scriptsize\rm I}\times\mbox{\bf
1}_{\mbox\scriptsize\rm O})+ \mbox{\bf 1}_{\mbox\scriptsize\rm
I}$. Then for $f^n$ being

\[
(l^1_{A,{\mbox\scriptsize\rm I}}\times\mbox{\bf
1}_{\mbox\scriptsize\rm O})^n_{\mbox\scriptsize\rm I}\cirk
k^1_{(A\times {\mbox\scriptsize\rm O})^n_{\mbox\tiny\rm I},
{\mbox\scriptsize\rm O}}: A^{n+1}_{\mbox\scriptsize\rm O}\vdash
A^{n+1}_{\mbox\scriptsize\rm I}
\]

\noindent and $g^n$ being

\[
l^1_{(A+ {\mbox\scriptsize\rm I})^n_{\mbox\tiny\rm O},
{\mbox\scriptsize\rm I}}\cirk (k^1_{A,{\mbox\scriptsize\rm
O}}+\mbox{\bf 1}_{\mbox\scriptsize\rm I})^n_{\mbox\scriptsize\rm
O}: A^{n+1}_{\mbox\scriptsize\rm O}\vdash
A^{n+1}_{\mbox\scriptsize\rm I}
\]

\noindent we have $G(f^n)=G(g^n)$, but we suppose that $f^n=g^n$
does not hold in $\cal D$. The equation $f^0=g^0$ is

\begin{tabbing}

\quad$((l^1_{A,{\mbox\scriptsize\rm I}}\times\mbox{\bf
1}_{\mbox\scriptsize\rm O})+\mbox{\bf 1}_{\mbox\scriptsize\rm
I})\cirk k^1_{(A\times {\mbox\scriptsize\rm O})+{\mbox\tiny\rm I},
{\mbox\scriptsize\rm O}}=l^1_{(A+ {\mbox\scriptsize\rm
I})\times{\mbox\tiny\rm O}, {\mbox\scriptsize\rm I}}\cirk
(k^1_{A,{\mbox\scriptsize\rm O}}+\mbox{\bf 1}_{\mbox\scriptsize\rm
I})\times\mbox{\bf 1}_{\mbox\scriptsize\rm O}:$
\\[.5ex]
\`$((A\times\mbox{\rm O})+\mbox{\rm I})\times\mbox{\rm O}\vdash
((A+\mbox{\rm I})\times\mbox{\rm O})+\mbox{\rm I}.$
\end{tabbing}

\noindent Note that $A^{n+1}_{\mbox\scriptsize\rm O}$ is not
O-normal, and  $A^{n+1}_{\mbox\scriptsize\rm I}$ is not I-normal.

We don't know whether it is sufficient to add to $\cal D$ the
equations $f^n=g^n$ for every $n\geq 0$ in order to obtain full
coherence for the resulting category.

\end{document}